\documentclass{article}
\usepackage{amssymb,epsf,graphicx}
\def\msy{\mathbb}
\newcommand{\qed}{\unskip\nobreak\hfill\mbox{ $\Box$}\par}
\def\bbbr{{\msy R}}
\def\bbbc{{\msy C}}
\def\bbbq{{\msy Q}}

\def\bbbz{{\msy Z}}
\def\bbbf{{\msy F}}

\def\bbbfa{\overline{\bbbf}}
\def\v#1{{\bf#1}}
\def\is{\equiv}
\def\mod#1{({\rm mod}\ #1)}

\newtheorem{theorem}[subsection]{Theorem}
\newtheorem{definition}[subsection]{Definition}

\newtheorem{proposition}[subsection]{Proposition}

\def\v#1{{\bf #1}}

\begin{document}

\title{Algebraic A-hypergeometric functions}
\author{Frits Beukers}
\date{December 5, 2008}
\maketitle

\abstract{We formulate and prove a combinatorial criterion to decide
if an A-hypergeometric system of differential equations has a full
set of algebraic solutions or not. This criterion generalises the so-called
interlacing criterion in the case of hypergeometric functions of one variable.}

\parindent=0pt

\section{Introduction}
The classically known hypergeometric functions of Euler-Gauss ($_2F_1$), its one-variable
generalisations $_{p+1}F_p$ and the many variable generalisations, such as
Appell's functions, the Lauricella functions and Horn series are all examples
of the so-called $A$-hypergeometric functions introduced by Gel'fand, Kapranov,
Zelevinsky in \cite{1,2,3}. We like to add that completely independently
B.Dwork developed a theory of generalised hypergeometric functions in \cite{13}
which is in many aspects parallel to the theory of A-hypergeometric functions.
The connection between the theories has been investigated in \cite{5} and \cite{14}.

The definition of A-hypergeometric functions begins with a finite subset
$A\subset\bbbz^r$ (hence their name)
consisting of
$N$ vectors $\v a_1,\ldots,\v a_N$ such that
\begin{itemize}
\item[i)] The $\bbbz$-span of $\v a_1,\ldots,\v a_N$ equals $\bbbz^r$.
\item[ii)]
There exists a linear form $h$ on $\bbbr^r$ such that $h(\v a_i)=1$ for all $i$.
\end{itemize}
The second condition ensures that we shall be working in the case of so-called
Fuchsian systems. In a number of papers, eg \cite{5}, this condition is dropped to
include the case of so-called confluent hypergeometric equations. 

We are also given a vector of parameters
$\alpha=(\alpha_1,\ldots,\alpha_r)$ which could be chosen in $\bbbc^r$, but we shall
restrict to $\alpha\in\bbbr^r$. The lattice $L\subset\bbbz^N$ of relations
consists of all $(l_1,\ldots,l_N)\in\bbbz^N$ such that $\sum_{i=1}^Nl_i\v a_i=0$.

The A-hypergeometric equations are a set of partial differential equations with
independent variables $v_1,\ldots,v_N$. This set consists of two groups. The first
are the structure equations

$$\Box_{\v l}\Phi:=\prod_{l_i>0}\partial_i^{l_i}\Phi-\prod_{l_i<0}\partial_i^{|l_i|}\Phi=0\eqno{{\rm (A1)}}$$
for all $\v l=(l_1,\ldots,l_N)\in L$.

The operators $\Box_{\v l}$ are called the box-operators.
The second group consists of the homogeneity or Euler equations.

$$Z_i\Phi:=\left(a_{1,i}v_1\partial_1+a_{2,i}v_2\partial_2+\cdots
+a_{N,i}v_N\partial_N-\alpha_i\right)\Phi =0,\ i=1,2,\ldots,r \eqno{{\rm
(A2)}}$$ where $a_{k,i}$ denotes the $i$-th coordinate of $\v a_k$.

In general the A-hypergeometric system is a holonomic system of dimension equal to the $r-1$-dimensional volume
of the so-called A-polytope $Q(A)$. This polytope is the convex hull of the endpoints of the $\v a_i$.
The volume-measure is normalised to $1$ for a $r-1$-simplex of lattice-points in the plane $h(\v x)=1$
having no other lattice points in its interior. In the first days of the theory of
A-hypergeometric systems there was some confusion as to what 'general' means, see \cite{5}. To avoid these
difficulties we make an additional assumption, which ensures that the dimension of the A-hypergeometric system
indeed equals the volume of $Q(A)$

\begin{itemize}
\item[iii)] The $\bbbr_{\ge0}$-span of $A$ intersected with $\bbbz^r$ equals the
$\bbbz_{\ge0}$-span of $A$.
\end{itemize}

Under this condition we have the following Theorem.

\begin{theorem}[GKZ, Adolphson]
Let notations be as above. If condition {\rm (iii)} is satisfied then the system of A-hypergeometric
differential equations is holonomic of rank equal to the volume of the convex hull $Q(A)$ of
$A$.
\end{theorem}

For a complete story on the dimension of the solution space we refer to \cite{18}.
In the present paper we shall use condition (iii) in the proof of the important
Proposition \ref{modpSolutions1}. 

To describe the standard hypergeometric solution of the A-hypergeometric system we define the projection
map $\psi_L:\bbbr^N\to\bbbr^r$ given by $\psi_L:\v e_i\mapsto\v a_i$ for $i=1,\ldots,N$.
Here $\v e_i$ denotes the $i$-th vector in the standardbasis of $\bbbr^N$.
Clearly the kernel of $\psi_L$ is the space $L\otimes\bbbr$. Choose a point
$\gamma=(\gamma_1,\ldots,\gamma_N)$ in $\psi_L^{-1}(\alpha)$, in other words choose $\gamma_1,
\ldots,\gamma_N$ such that $\alpha=\gamma_1\v a_1+\cdots+\gamma_N\v a_N$. Then a formal solution
of the A-hypergeometric system can be given by
$$\Phi_{L,\gamma}(v_1,\ldots,v_N)=\sum_{\v l\in L}{\v v^{\v l}\over\Gamma(\v l+\gamma+\v 1)}$$
where we use the short-hand notation
$${\v v^{\v l+\gamma}\over\Gamma(\v l+\gamma+\v 1)}={v_1^{l_1+\gamma_1}\cdots v_N^{l_N+\gamma_N}
\over \Gamma(l_1+\gamma_1+1)\cdots\Gamma(l_N+\gamma_N+1)}.$$
By a proper choice of $\gamma\in\psi_L^{-1}(\alpha)$ this formal solution gives rise
to actual powerseries solutions with a non-trivial region of convergence.

The real positive cone generated by the vectors $\v a_i$ is denoted by $C(A)$. This is a
polyhedral cone with a finite number of faces. We recall the following
Theorem.

\begin{theorem}\label{irreducible} The A-hypergeometric system is irreducible if and only if
$\alpha+\bbbz^r$ contains no points in any face of $C(A)$.
\end{theorem}

This Theorem is proved in \cite{4}, Theorem 2.11 using perverse sheaves. It would be nice to have a
more elementary proof however.

Let us now assume that $\alpha\in\bbbq^r$. We shall be interested
in those irreducible A-hypergeometric system that have a complete set of
solutions algebraic over $\bbbc(v_1,\ldots,v_N)$. This question
was first raised in the case of Euler-Gauss hypergeometric functions and
the answer is provided by the famous list of H.A.Schwarz, see
\cite{6}. In 1989 this list was extended to general one-variable
$_{p+1}F_p$ by Beukers and Heckman, see \cite{7}. For the several
variable cases, a characterization for Appell-Lauricella $F_D$ was
provided by Sasaki \cite{8} in 1976 and Wolfart, Cohen \cite{9} in 1992. The
Appell systems $F_4$ and $F_2$ were classified by M.Kato in \cite{10} (1997)
and \cite{11} (2000).

In the case of one-variable hypergeometric functions there is a simple combinatorial
criterion to decide if they are algebraic or not. Consider the hypergeometric function
$$_pF_{p-1}(\alpha_1,\ldots,\alpha_p;\beta_1,\ldots,\beta_{p-1}|z).$$
Define $\beta_p=1$. We assume $\alpha_i-\beta_j\not\in\bbbz$ for all $i,j=1,\ldots,p$,
which ensures that the corresponding hypergeometric differential equation is irreducible.
We shall say that the sets $\alpha_i$ and $\beta_j$ interlace modulo 1 if the points
of the sets $e^{2\pi i\alpha_j}$ and $e^{2\pi i\beta_k}$ occur alternatingly when
running along the unit circle. The following Theorem is proved in \cite{7}.

\begin{theorem}[Beukers, Heckman] Suppose the one-variable hypergeometric equation with
parameters $\alpha_i,\beta_j\in\bbbq$ ($i,j=1,2,\ldots,p$) with $\beta_p=1$ is irreducible.
Let $D$ be the common denominator of the parameters. Then the solution set of the
hypergeometric equation consists of algebraic functions (over $\bbbc(z)$) if and only
if the sets $k\alpha_i$ and $k\beta_j$ interlace modulo 1 for every integer $k$
with $1\le k<D$ and $\gcd(k,D)=1$.
\end{theorem}

It is the purpose of this paper to generalize the interlacing
condition to a similar condition for A-hypergeometric systems. We assume that
$\alpha\in\bbbr^r$ and define
$K_{\alpha}=(\alpha+\bbbz^r)\cap C(A)$. A point $\v p\in
K_{\alpha}$ is called an {\it apexpoint} if  $\v p\not\in\v
q+C(A)$ for every $\v q\in K_{\alpha}$ with $\v q\ne\v p$. We call
the number of apexpoints the {\it signature} of the polytope
$A$ and parameters $\alpha$. Notation: $\sigma(A,\alpha)$. 

\begin{proposition}\label{apexbound} Let $\alpha\in\bbbr^r$. Then 
$\sigma(A,\alpha)$ is less than or equal to the volume
of the $A$-polytope $Q(A)$.
\end{proposition}

We say that the signature is {\it maximal} if it equals the volume of $Q(A)$.

\begin{theorem}\label{algebraic} Let $\alpha\in\bbbq^r$ and suppose the A-hypergeometric system is irreducible.
Let $D$ be the common denominator of the coordinates of $\alpha$. Then the solution
set of the
A-hypergeometric system consists of algebraic solutions (over $\bbbc(v_1,\ldots,v_N)$) if and
only if $\sigma(A,k\alpha)$ is maximal
for all integers $k$ with $1\le k<D$ and $\gcd(k,D)=1$.
\end{theorem}

To compare this result with the one-variable interlacing condition for $_2F_1$
we illustrate a connection. In the case of Euler-Gauss hypergeometric function we
have $r=3,N=4$ and
$$
\v a_1=(1,0,0),\ \v a_2=(0,1,0),\ \v a_3=(0,0,1),\ \v a_4=(1,1,-1).
$$
The faces of the cone
generated by $\v a_i$ ($i=1,\ldots,4$) are given by $x=0,y=0,x+z=0,y+z=0$
(we use the coordinates $x,y,z$ in $\bbbr^3$. We define
$\alpha=(-a,-b,c-1)$. Theorem \ref{irreducible} implies that
irreducibility comes down to the inequalities
$-a,-b,-a+c,-b+c\not\in\bbbz$. These are the familiar irreducibility conditions
for the Euler-Gauss hypergeometric functions.

The lattice of relations has rank one and is generated by $(-1,-1,1,1)$.
We choose $\gamma=(-a,-b,c-1,0)$. Then the
formal solution $\Phi_{L,\gamma}$ reads
$$v_1^{-a} v_2^{-b} v_3^{c-1}\sum_{k\in\bbbz}
{v_1^{-k} v_2^{-k} v_3^k v_4^k\over\Gamma(-k-a+1)
\Gamma(-k-b+1)\Gamma(c+k)\Gamma(k+1)}.$$
Clearly $1/\Gamma(k+1)$ vanishes for $k\in\bbbz_{<0}$, so our
summation actually runs over $k\in\bbbz_{\ge0}$. Apply the
identity $1/\Gamma(1-z)=\sin(\pi z)\Gamma(z)/\pi$ to $z=k+a$
and $z=k+b$ to obtain
$$\Phi_{L,\gamma}=v_1^{-a} v_2^{-b} v_3^{c-1}{\sin(\pi a)\sin(\pi b)\over\pi^2}
\sum_{k=0}^{\infty}{\Gamma(a+k)\Gamma(b+k)\over\Gamma(c+k)k!}\left(
{v_3v_4\over v_1v_2}\right)^k.$$
Setting $v_1=v_2=v_3=1$ and $v_4=z$ we recognize the Euler-Gaussian hypergeometric
series $_2F_1(a,b,c|z)$. By shifting over $\bbbz$ if necessary we can see to it
that $a,b,c$ are in the interval $(-1,0)$. Suppose that the sets $\{a,b\}$
and $\{0,c\}$ interlace modulo $1$. By interchange of $a,b$ if necessary
we can restrict ourselves to the case $-1<a<c<b<0$. It is straightforward
to verify that $(-a,1-b,c)$ and $(-a,-b,1+c)$ are apexpoints of $K_{\alpha}$.
If the sets do not interlace then one checks that $(-a,-b,c)$ is the unique
apexpoint if $a,b<c$ and $(-a,-b,1+c)$ is the unique apexpoint if $c<a,b$.

A second example is Appell's hypergeometric equation $F_2$. The Appell $F_2$
hypergeometric function reads
$$F_2(a,b,b',c,c'|x,y)=\sum_{m=0}^{\infty}\sum_{n=0}^{\infty}
{(a)_{m+n}(b)_m(b')_n\over (c)_m(c')_n m! n!}\ x^m y^n,$$
where $(x)_n$ denotes the Pochhammer symbol defined by 
$\Gamma(x+n)/\Gamma(x)=x(x+1)\cdots(x+n-1)$.
The function $F_2$ satisfies a system of partial differential equations
of rank 4.
Algebraicity of these functions is completely described in \cite{11}.

The A-parameters are as follows. We have $r=5$ and $N=7$.
The set $A$ consists of the standard basisvectors $\v a_1,\ldots,\v a_5$
in $\bbbr^5$ and $\v a_6=(1,1,0,-1,0),\ \v a_7=(1,0,1,0,-1)$.
We take $\alpha=(-a,-b,-b',c,c')$. The lattice
$L$ of relations is generated by $(-1,-1,0,1,0,1,0)$ and $(-1,0,-1,0,1,0,1)$.
Take $\gamma=(-a,-b,-b',c,c',0,0)$.
In a similar way as we did for the Euler-Gauss functions we can now go from the
formal expansion $\Phi_{L,\gamma}$ to the explicit Appell function $F_2$.

One can compute that the cone $C(A)$ has 8 faces and they
are given by $x_1=0, x_2=0, x_3=0$, $x_2+x_4=0, x_3+x_5=0, x_1+x_4=0,
x_1+x_5=0$ and $x_1+x_4+x_5=0$. Using Theorem \ref{irreducible} it follows
that the A-hypergeometric system is irreducible if and only if none of the following numbers
is an integer,
$$a,b,b',-b+c,-b'+c',-a+c,-a+c',-a+c+c'.$$
These are precisely the irreducibility condtions for the $F_2$-system given
in \cite{11}. In that paper it is shown for example that with the
choice
$$\alpha=(-a,-b,-b',c,c')=(1/10, 7/10, 9/10, 3/5, 1/5)$$
the solutions of the rank 4 Appell system are all algebraic. A small computer
calculation shows that the apex points of $K_{\alpha}$ are given by
$$(21/10, 7/10, 9/10, -2/5,-4/5),\quad (11/10, 7/10, 9/10, 3/5, -4/5)$$
$$11/10, 7/10, 9/10, -2/5, 1/5),\quad (1/10, 7/10, 9/10, 3/5, 1/5).$$
Similarly there are four apexpoints for the conjugate parameter 5-tuples
$3\alpha,7\alpha,9\alpha$.

\section{A simple example}
In this section we show a more elaborate example of algebraic hypergeometric
functions of Horn-type G3 which, to our knowledge, has not been dealt with
before. The corresponding series with parameters $a,b$ is given by
$$G_3(a,b,x,y)=\sum_{m\ge0,n\ge0}{(a)_{2m-n}(b)_{2n-m}\over m!n!}x^my^n.$$
Here again, $(x)_n$ denotes the Pochhammer symbol defined by
$\Gamma(x+n)/\Gamma(x)$. However, now the index $n$ may
be negative, in which case the definition explicitly reads 
$\Gamma(x+n)/\Gamma(x)=1/(x-1)\cdots(x-|n|)$. The system of
differential equations is a rank 3 system.
The set $A$ can be chosen in $\bbbz^2$ for example as
$$\v a_1=(-1,2),\quad \v a_2=(0,1),\quad \v a_3=(1,0),\qquad \v a_4=(2,-1).$$
Below we show a picture of the cone $C(A)$ spanned by the elements of $A$,
together with the points from $A$. In addition, the dark grey area indicates
the set of apexpoints with respect to $A$. The parameter vector of the
corresponding A-hypergeometric system is given by $(-a,-b)$. 

\begin{theorem}\label{G3}
Consider the A-hypergeometric system corresponding to the Horn G3 equations. In the
following cases the system is irreducible and has only algebraic solutions.
\begin{enumerate}
\item $a+b\in\bbbz$ and $a,b\not\in\bbbz$.
\item $a\is1/2\mod{\bbbz}, b\is 1/3,2/3\mod{\bbbz}$
\item $a\is1/3,2/3\mod{\bbbz}, b\is1/2\mod{\bbbz}$
\end{enumerate}
\end{theorem}

{\bf Proof} 
This is an application of Theorem \ref{algebraic}. In all cases we need
only be interested in $a,b\mod{\bbbz}$. 
In the following picture the light grey area is the cone $C(A)$,
the dark grey area indicates the location
of the apexpoints. If $a+b\in\bbbz$ then we note that 
there are precisely 3 points
of $(-a,-b)+\bbbz^2$ in the dark grey area, all lying on the line $x+y=1$.
Hence three apexpoints.
\medskip

\centerline{\includegraphics[height=4.5cm]{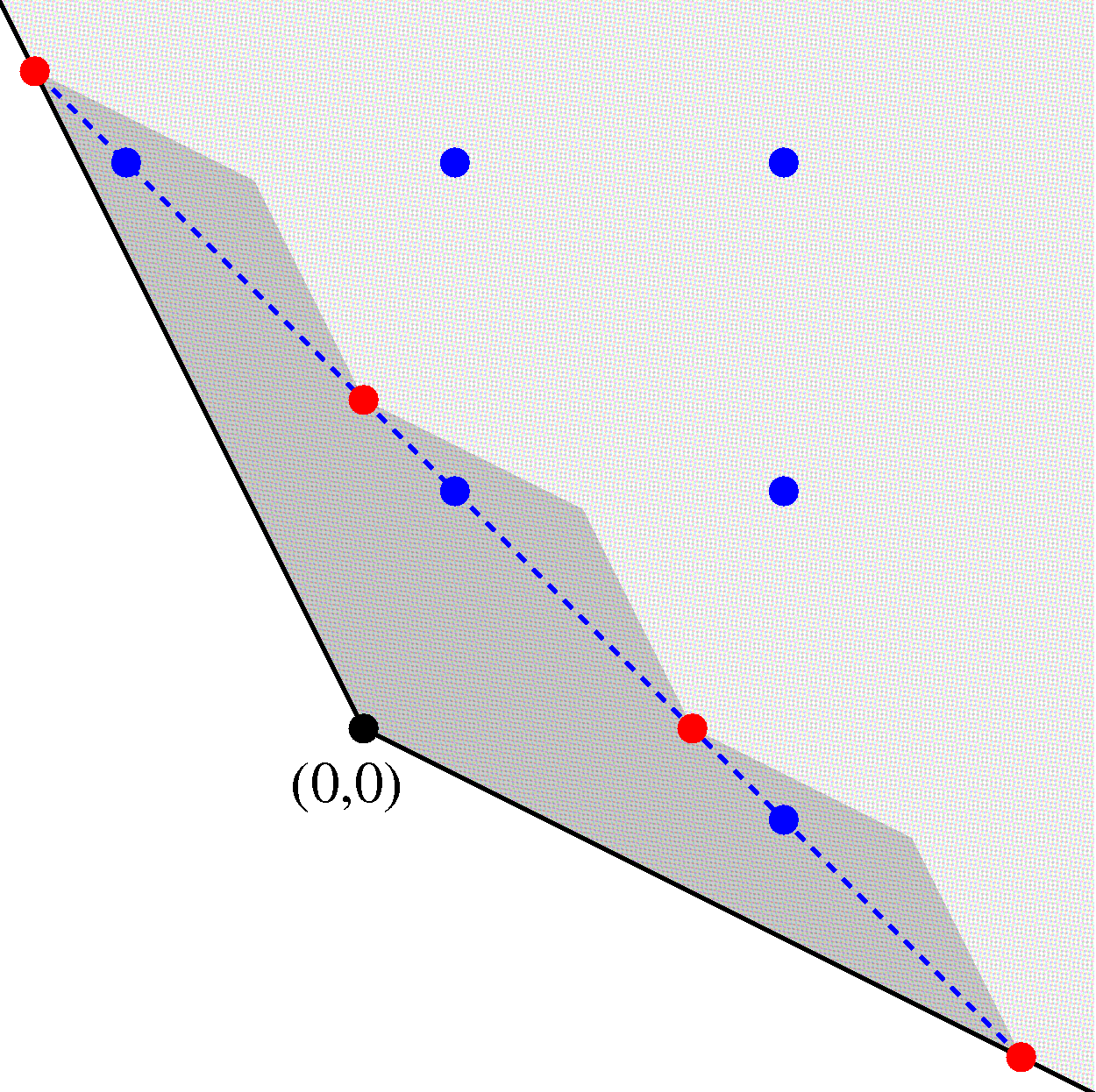}}

If $a+b\in\bbbz$ then also $ka+kb\in\bbbz$ for any integer. 
Irreducibility of the systems is ensured by Theorem \ref{irreducible}
and the fact that $a,b\not\in\bbbz$. Therefore, in
the first case all conditions of Theorem \ref{algebraic} are
fulfilled. 

In the following picture we have drawn the sets $(1/2,1/3)+\bbbz^2$
and $(1/2,2/3)+\bbbz^2$ intersected with $C(A)$.
\medskip

\centerline{\includegraphics[height=4.5cm]{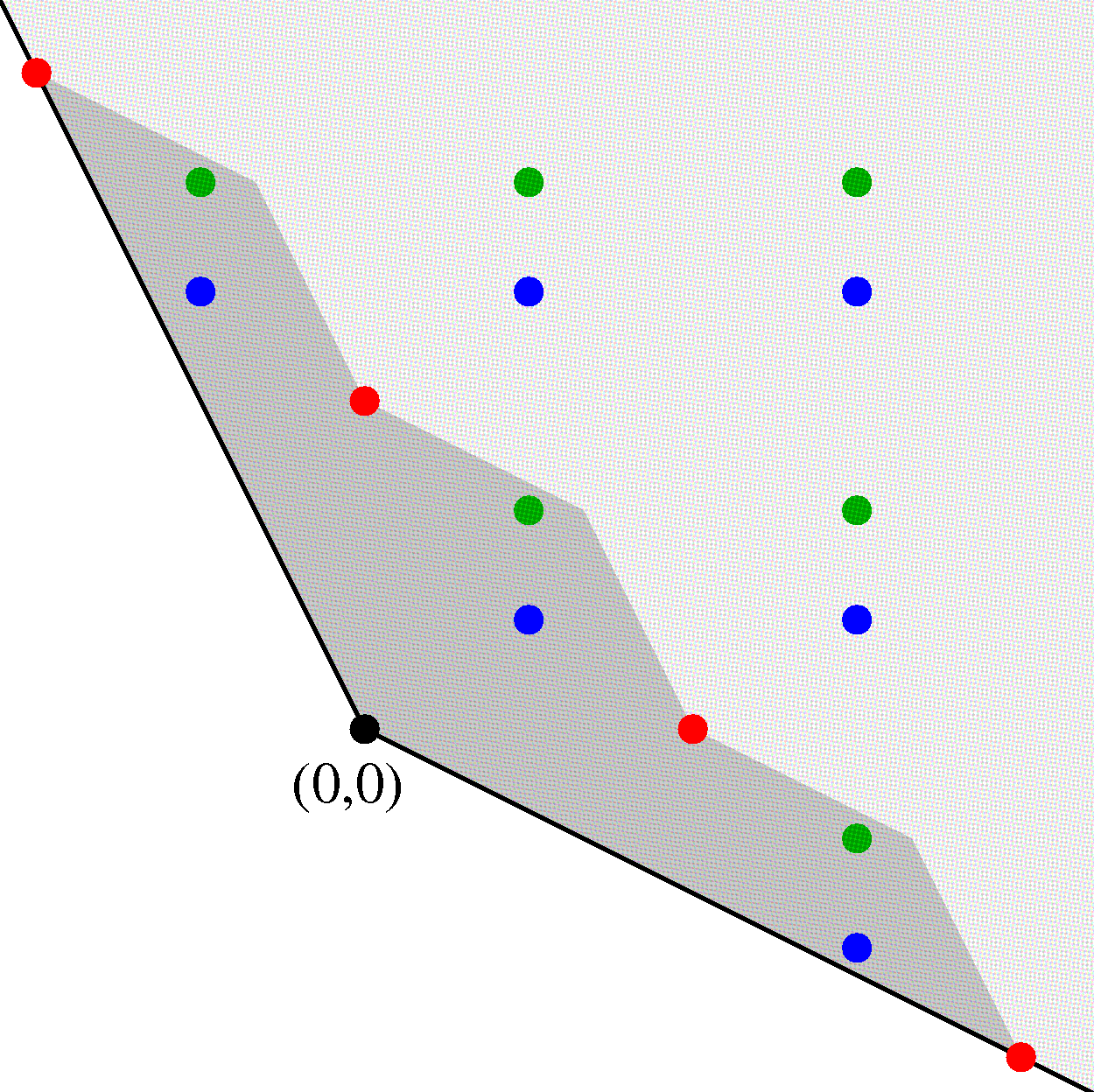}}

Clearly each set 
has three apexpoints and Theorem \ref{algebraic} can be applied to
prove the second case. The third case runs similarly.

\qed

\medskip

It has been verified by J.Schipper, an Utrecht graduate student, that 
Theorem \ref{G3} gives the characterisation of all irreducible Horn G3-systems
with algebraic solutions.

A fairly involved calculation reveals that a formula for
$G_3(a,1-a,x,y)$ can be given as follows

$$G_3(a,1-a,x,y)=f(x,y)^a\sqrt{{g(x,y)\over \Delta}}$$
where
$$\Delta=1+4x+4y+18xy-27x^2y^2$$
and 
$$xf^3-y=f-f^2,\qquad g(g-1-3x)^2=x^2\Delta.$$
For reference we display the series expansions of $f$ and $g$.
\begin{eqnarray*}
f&=&1+(y-x) +\left(2 x^2-y x-y^2\right) +\left(-5 x^3+3
   y x^2+2 y^3\right)\\
   &&+\left(14 x^4-10 y x^3+y^3 x-5y^4\right)+O(x,y)^5
   \end{eqnarray*}
$$g=1+2 x -x(x+2 y)+2x(x-y)^2 -x(x-y)^2(5 x+4
   y) +O(x,y)^5$$

Moreover, $g=1+4x-2xf-3x^2f^2$. In particular,
$f$ and $g$ generate the same cubic
extension of $\bbbq(x,y)$. 

\section{The signature}
{\bf Proof} of Proposition \ref{apexbound}. We use the following property.
Let $\v b_1,\ldots,\v b_r\in\bbbz^r$ be independent vectors. Let $\beta\in\bbbr^r$.
Then the number of points of $\beta+\bbbz^r$ inside the fundamental block
$\{\sum_{i=1}^r\lambda_i\v b_i\ |\ 0\le\lambda_i<1\}$ is equal to
$|\det(\v b_1,\ldots,\v b_r)|$.

Write $Q(A)$ as a union of
$r-1$-simplices $\cup_{i=1}^m\sigma_i$ (a so-called triangulation of $Q(A)$).
Every simplex $\sigma_i$ is spanned by $r$ independent vectors $\v a_{i_j}$
($j=1,\ldots,r$). Let $B_i$ be the fundamental block spanned by these vectors.

Let $\v a$ be an apexpoint of $K_{\alpha}$. Then $\v a$ is
contained in a positive cone spanned by one of the simplices
$\sigma_i$. For every choice of $\v a_{i_j}$ ($j=1,\ldots,r$) the
point $\v a-\v a_{i_j}$ falls outside this cone. If not, then $\v
a$ would be contained in $\v a_{i_j}+C(A)$. Hence $\v a\in B_i$.
Since $B_i$ contains at most $|\det(B_i)|$ point from
$\alpha+\bbbz^r$, we see that the number of apexpoints is bounded
above by $\sum_{i=1}^m|\det(B_i)|$. This equals precisely the
$r-1$-dimensional volume of the $A$-polytope $Q(A)$.

\qed
\medskip

For any $\v k=(k_1,\ldots,k_N)\in\bbbz^N$ we define the hypercube
$$F(\v k)=\{(x_1,\ldots,x_N)| k_i\le x_i<k_i+1,\ (i=1,2,\ldots,N)\}$$
The space $\bbbr^N$ can be seen as the disjoint union of cells
$F(\v k)$ with $\v k$ running over $\bbbz^N$. Let us denote
$L(\bbbr)=L\otimes\bbbr$. We intersect the union $\cup_{\v k}F(\v
k)$ with the translated space $\gamma+L(\bbbr)$. Each hypercube
intersects $\gamma+L(\bbbr)$ in a cell which may, or may not be
closed in $\gamma+L(\bbbr)$. Let us denote $V(\v k)= F(\v
k)\cap(\gamma+L)$. We call $V(\v k)$ a compact cell if it is
closed and non-empty. Of course, if we shift a compact cell over a
point of $L$, we get another compact cell. In the following
Proposition we denote the shifted hyperquadrant $\v
x+\bbbr_{\ge0}^N$ by $P(\v x)$.

\begin{proposition}\label{compactcell}
With the notations as above, $V(\v k)$ is a compact cell if and only if
$P(\v k)$ has non-trivial intersection with $\gamma+L$ and
$P(\v k+\v e_i)$ has empty intersection with $\gamma+L$
for $i=1,2,\ldots,N$. In particular, $V(\v k)=P(\v k)\cap(\gamma+L(\bbbr))$.
\end{proposition}

{\bf Proof}. Let us denote the intersection of $\cup_{i=1}^N P(\v k+\v e_i)$
with $\gamma+L$ by $W$. Notice that $W$ is the set-theoretic
difference between $P(\v k)\cap(\gamma+L)$ and $V(\v k)$. In particular
$W$ is a closed set. Suppose that $V(\v k)$ is a compact cell.
The only way that the difference $W$ of the two non-empty closed convex
sets $P(\v k)\cap(\gamma+L)$ and $V(\v k)$ can be closed is when $W$
is empty. Hence $P(\v k+\v e_i)\cap(\gamma+L)$ is empty for all $i$.

Suppose conversely that $W$ is empty. Then $V(\v k)=P(\v k)$. Since
$P(\v k)$ is closed, the same should hold for $V(\v k)$. Since $V(\v k)$
is also bounded, we conclude that $V(\v k)$ is compact.

\qed
\medskip

In the following recall the map $\psi:\bbbr^N\to\bbbr^r$ given by 
$\psi:\v e_i\mapsto \v a_i$ for $i=1,\ldots,N$.

\begin{proposition}\label{cell_apex}
The compact cells in $\gamma+L(\bbbr)$, modulo $L$,
are in 1-1 correspondence with the apex-points of $K_{\alpha}$.
The correspondence is given by $V(\v k)\mapsto\alpha-\psi(\v k)$.

Let $\v a,\v a'$ be two different apexpoints. Then $\psi^{-1}(\v
a)\cap P(\v 0)$ and $\psi^{-1}(\v a')\cap P(\v 0)$ are disjoint
and contained in the unit cube $0\le x_i<1$ for $i=1,2,\ldots,N$.

\end{proposition}

{\bf Proof} An apexpoint $\v a\in\alpha+\bbbz^r$ is characterized
by the fact that $\v a\in C(A)$ and $\v a-\v a_i\not\in C(A)$ for
$i=1,\ldots, N$. Notice that $\psi:\bbbr^N\to\bbbr^r$ is actually the quotient map
$\bbbr^N\to\bbbr^N/L(\bbbr)$. Let $\v k\in\bbbz^N$. First recall that $\psi(\gamma)=\alpha$ and
observe that $\psi(P(\v k))=\psi(\v k)+C(A)$.
As a result we see that $P(\v k)\cap(\gamma+L(\bbbr))$ is non-empty if
and only if $\alpha-\psi(\v k)\in C(A)$. Let us assume that $V(\v k)$ is
compact and apply Proposition
\ref{compactcell}. Then
$\alpha-\psi(\v k)\in C(A)$ and $\alpha-\psi(\v k)-\v a_i\not\in
C(A)$ for $i=1,\ldots,N$. Hence $\alpha-\psi(\v k)$ is an
apexpoint. 

Conversely, when $\v a$ is an apexpoint, find $\v
k\in\bbbz^N$ such that $\v a=\alpha-\psi(\v k)$. Then we find that
$(\gamma+L(\bbbr))\cap P(\v k)$ is non empty and the
sets $\psi^{-1}(\alpha)\cap P(\v k+\v e_i)$ for $i=1,2,\ldots,N$
are empty. Therefore, by application of Proposition \ref{compactcell},
$V(\v k)$ is a compact cell.

Let $\v a$ be an apexpoint and choose $\v k\in\bbbz^N$ such that
$\alpha-\psi(\v k)=\v a$. Then $\psi^{-1}(\alpha)\cap P(\v k)$ is
a compact cell. Consequently, after shifting over $\v k$, the set
$\psi^{-1}(\v a)\cap P(\v 0)$ is a compact cell in $\gamma-\v
k+L(\bbbr)$. Hence $\psi^{-1}(\v a)\cap P(\v 0)$ is contained in the unit cube
in $\bbbz^N$. Two sets $\psi^{-1}(\v a)$ and $\psi^{-1}(\v a')$
are obviously distinct whenever $\v a$ and $\v a'$ are distinct.

\qed
\medskip

\section{Mod $p$ solutions}
Let us assume that $\alpha\in\bbbz^r$ and let $p$ be a prime. We
describe the polynomial solutions in $\bbbf_p[v_1,\ldots,v_N]$ of
the A-hypergeometric system with parameters $A$ and $\alpha$ considered modulo
$p$.

Let $\beta_i/p$ ($i=1,\ldots,\sigma$) be the set of apexpoints of
$(\alpha/p+\bbbz^r)\cap C(A)$. To any apexpoint we associate the
set of lattice points
$$\Gamma_i=\psi^{-1}(\beta_i)\cap\bbbz_{\ge0}^N$$
where $\psi:\bbbr^N\to \bbbr^N/L(\bbbr)$ is defined as in the previous section.
For any $i=1,2,\ldots,\sigma$ we define
$$\Psi_i:=\sum_{\v l\in\Gamma_i}{{\v v}^{\v l}\over\Gamma(\v l+1)}.$$
This is a polynomial solution to the A-hypergeometric system with parameters
$A,\beta_i$. Since $\beta_i/p$ is an apexpoint of
$(\alpha/p+\bbbz^r)\cap C(A)$, the preimage
$\psi^{-1}(\beta_i/p)\cap P(\v 0)$ is contained in the unit cube
in $\bbbr^N$ according to Proposition \ref{cell_apex}. Hence the
points of $\Gamma_i=\psi^{-1}(\beta_i)\cap\bbbz_{\ge0}^N$ are
contained in the cube $0\le x_i<p$ for $i=1,\ldots,N$. In
particular none of the positive coordinates of any $\v
l\in\Gamma_i$ is divisible by $p$, hence $\Gamma(\v l+\v
1)\not\is0\mod{p}$ for all $\v l\in\Gamma_i$. This means that
$\Psi_i$ can be reduced modulo $p$. Furthermore, since
$\alpha\is\beta_i\mod{p}$, the polynomial $\Psi_i$ is a
polynomial solution modulo $p$ to the A-hypergeometric system with parameters
$A,\alpha$.

Since each of the sets $\Gamma_i$ is contained in the cube $0\le
x_i<p$ any two shifts $\Gamma_i+p\v k_i$ and $\Gamma_j+p\v k_j$
for different $i,j$ and $\v k_i,\v k_j\in\bbbz^N$ are disjoint. In
particular the polynomials $\Psi_i$ are independent over
$\bbbf_p[v_1^p,\ldots,v_N^p]$.

\begin{proposition}\label{modpSolutions1}
Every mod $p$ polynomial solution of the A-hypergeometric system with
parameters $A,\alpha\in\bbbz^r$ is an
$\bbbf_p[v_1^p,\ldots,v_N^p]$-linear combination of the
polynomials $\Psi_i$.
\end{proposition}

Let $P=\sum_{\v m}p_{\v m}\v v^{\v m}$ be a polynomial solution of
the A-hypergeometric system with parameters $A,\alpha$. For formal reasons we
extend the summation over all of $\bbbz^N$ but it should be
understood that the set of multi-indices $\v m$ with $p_{\v
m}\not\is0\mod{p}$ is finite and contained in $\bbbz_{\ge0}^N$.

Let us substitute this in the system (A1). Any $\v l\in L$ ca be decomposed
as $\v l=\v l_+-\v l_-$ where $\v l_+,\v l_-\in\bbbz_{\ge0}^N$ and we assume
they have disjoint support. Denoting $[\v m]_{\v
r}=\prod_{i=1}^N m_i(m_i-1)\cdots(m_i-r_i+1)$ for any $\v
r\in\bbbz_{\ge0}^N$ we obtain
\begin{eqnarray*}
0&\is&\sum_{\v m}\left([\v m]_{l_+}p_{\v m}\v v^{\v m-l_+}-
[\v m]_{l_-}p_{\v m}\v v^{\v m-l_-}\right)\mod{p}\\
&\is&\v v^{-l_+}\sum_{\v m}\left([\v m]_{l_+}p_{\v m}\v -[\v m-\v l]_{l_-}p_{\v m-\v l}
\right)v^{\v m}\mod{p}
\end{eqnarray*}
Hence
\begin{equation}\label{recursion}
[\v m]_{l_+}p_{\v m} -[\v m-\v l]_{l_-}p_{\v m-\v l}\is0\mod{p}
\end{equation}
for every $\v m\in \bbbz^N$ and every $\v l\in L$. Substitution in
(A2) gives
$$\sum_{\v m}(-\alpha_i+a_{1i}m_1+\cdots+a_{Ni}m_N)\v p_{\v m}v^{\v m}\is0\mod{p}$$
for $i=1,2,\ldots,r$.

Hence $-\alpha_i+a_{1i}m_1+\cdots+a_{Ni}m_N\is0\mod{p}$
for every $\v m$ with $p_{\v m}\not\is0\mod{p}$. Note that the
system (A2) gives no extra relations between different $p_{\v
m}$. They only require that $\psi(\v m)\is\alpha\mod{p}$.

The system (A1) relates only those coefficients $p_{\v m}$ for
which the multi-indexes $\v m$ differ by an element of $L$. Hence
we can split $P$ as a sum of terms of the form
$P_{\beta}=\sum_{\psi(\v m)=\beta}p_{\v m}\v v^{\v m}$ with
$\beta\in\bbbz^r$ and each summand $P_{\beta}$ satisfies modulo
$p$ the A-hypergeometric system with parameters $A,\alpha$. The equations
(A2) applied to $P_{\beta}$ tell us that
$\beta\is\alpha\mod{p}$.

We define a partial ordering on $\bbbr^N$. We say that $\v y\ge \v
x$ if all components of $\v y$ are larger or equal than the
corresponding component of $\v x$. In particular, when $\v y\ge \v
x$ and $\v y\ne\v x$ we write $\v y>\v x$.

When $\v m=(m_1,\ldots,m_N)\in\bbbz^N$ we denote by $\lfloor \v
m/p\rfloor$ the vector $(\lfloor m_1/p\rfloor,\ldots,\lfloor
m_N/p\rfloor)$. Consider the recursion (\ref{recursion}). We claim
that $[\v m]_{l_+}\is0\mod{p}$ if and only if $\lfloor \v
m/p\rfloor-\lfloor (\v m-\v l)/p\rfloor$ has at least one positive
component. This can be seen through the following sequence of
equivalences,

\begin{eqnarray*}
[\v m]_{\v l_+}\is0\mod{p}
&\iff& \exists i:\ m_i(m_i-1)\cdots(m_i-l_i+1)\is0\mod{p}\\
&\iff& \exists i,\lambda:\ 0\le\lambda<l_i,\ m_i-\lambda\is0\mod{p}\\
&\iff& \exists i,\lambda:\ 0\le\lambda<l_i,\ (m_i-\lambda)/p\in\bbbz\\
&\iff& \exists i:\ \lfloor m_i/p\rfloor -\lfloor(m_i-l_i)/p\rfloor>0
\end{eqnarray*}

Similarly we see that $[\v m-\v l]_{l_-}\is0\mod{p}$ if and
only if $\lfloor \v m/p\rfloor-\lfloor (\v m-\v l)/p\rfloor$ has
at least one negative component.

In terms of our partial ordering this implies that $[\v m]_{l_+}\is0\mod{p}$
and $[\v m-\v l]_{l_-}\is0\mod{p}$ if and only if neither
$\lfloor \v m/p\rfloor\ge \lfloor (\v m-\v l)/p\rfloor$ nor
$\lfloor \v m/p\rfloor\le \lfloor (\v m-\v l)/p\rfloor$, i.e
$\lfloor \v m/p\rfloor$ and $\lfloor (\v m-\v l)/p\rfloor$ are unrelated. 
Write $\v m'=\v m-\v l$, then $p_{\v m}$ and $p_{\v m'}$ are related through
(\ref{recursion}) if and only if $\psi(\v m)=\psi(\v m')$ and 
$\lfloor \v m/p\rfloor$ and $\lfloor \v m'/p\rfloor$
are related.

Now suppose that $p_{\v m}\not\is0\mod{p}$. We assert that for any
$\lambda\in L(\bbbr)$ the inequality $\lfloor \v m/p\rfloor \le
\lfloor (\v m-\lambda)/p\rfloor$ implies equality. First we deal
with the case when $\lambda=\v l\in L$. Suppose $\lfloor \v
m/p\rfloor < \lfloor (\v m-\v l)/p\rfloor$. Then $[\v
m]_{l_+}\not\is 0\mod{p}$ and $[\v m-\v l]_{l_-}\is0\mod{p}$. This
gives a contradiction with relation (\ref{recursion}). Hence
$\lfloor \v m/p\rfloor \le \lfloor (\v m-\v l)/p\rfloor$ implies
equality.

Now, in general, suppose that there exists $\lambda\in L(\bbbr)$
such that $\lfloor \v m/p\rfloor < \lfloor (\v
m-\lambda)/p\rfloor$. The vector $\v m-\lambda-p\lfloor(\v
m-\lambda)/p\rfloor$ has non-negative coefficients. Hence its
image under $\psi$ is contained in the cone $C(A)$. Moreover,
since $\psi(\lambda)=0$, the image has integer coordinates. Choose
a vector $\v k\in\bbbz_{\ge0}^N$ such that $\psi(\v k)=\psi(\v m-
p\lfloor(\v m-\lambda)/p\rfloor)$. Notice that this is only possible
because of Assumption iii) which we made in the introduction.
Hence there exists $\v l\in L$
such that $\v k=\v m-\v l-p\lfloor(\v m-\lambda)/p\rfloor$. In
particular, $\lfloor(\v m-\v l)/p\rfloor \ge \lfloor(\v
m-\lambda)/p\rfloor$. Since, by assumption, the latter vector is
strictly larger than $\lfloor\v m/p\rfloor$ we again get a
contradiction. Hence we conclude that
\begin{equation}\label{equality}
\lfloor(\v m-\lambda)/p\rfloor\ge \lfloor\v m/p\rfloor
\Rightarrow \lfloor(\v m-\lambda)/p\rfloor= \lfloor\v m/p\rfloor.
\end{equation}
Another
way of phrasing property (\ref{equality}) is to say that $\v m/p$
is contained in a compact cell of the affine space
$\v m/p+L(\bbbr)$. To see this consider the cell $V(\lfloor m/p\rfloor)$.
Of course it contains $\v m/p$. Let now $(\v m-\lambda)/p$ be any other
point in $P(\lfloor \v m/p\rfloor)\cap \v m/p+L(\bbbr)$. Then 
$\lfloor(\v m-\lambda)/p\rfloor \ge\lfloor\v m/p\rfloor$ and we have seen
that this implies equality. Therefore $(\v m-\lambda)/p$ is contained in
$V(\lfloor\v m/p\rfloor)$. Hence the latter cell is compact by Proposition
\ref{compactcell}. 

Let $\beta$ be as in the polynomial $P_{\beta}$ above and $\gamma\in\bbbr^N$
such that $\psi(\gamma)=\beta$. Let $\v
m,\v m'\in\bbbz_{\ge0}^N$ be such that $\psi(\v m)=\psi(\v
m')=\beta$ and such that $\v m/p,\v m'/p$ are in a compact cell of
$\gamma/p+L(\bbbr)$. Let $\v l=\v m-\v m'$. If $\v m/p,\v m'/p$
belong to different compact cells we have neither $\lfloor \v
m/p\rfloor \le \lfloor\v m'/p\rfloor$ nor $\lfloor \v m/p\rfloor
\ge \lfloor\v m'/p\rfloor$. Hence $p_{\v m}$ and $p_{\v m'}$ are
unrelated by relation (\ref{recursion}).

As a consequence of this all, the polynomial $P_{\beta}$ splits as
a sum of terms of the form $\sum_{\lfloor \v m/p\rfloor=\v k}p_{\v
m}\v v^{\v m}$ and each such sum is a solution of (A1) and
(A2). The latter summation can be rewritten as $\v v^{p\v
k}\sum_{\lfloor\v m/p\rfloor=\v k}p_{\v m}\v v^{\v m-p\v k}.$

The multi-indices $\v m$ in $\sum_{\lfloor\v m/p\rfloor=\v k}p_{\v
m}\v v^{\v m-p\v k}$ should in addition satisfy $\psi(\v
m)=\beta$. Replace $\v m$ by $\v n+p\v k$ and we obtain the
solution
\begin{equation}\label{reducedsolution}
\sum_{\lfloor\v n/p\rfloor=\v 0}b_{\v n}\v v^{\v n},
\end{equation}
where we put $b_{\v n}=p_{n+p\v k}$. We now know that all
multi-indices $\v n$ are contained in the cube $0\le x_i<p$ for
$i=1,\ldots,N$. Furthermore, in the recursion relation
$$
[\v n]_{l_+}b_{\v n} -[\v n-\v l]_{l_-}b_{\v n-\v l}\is0\mod{p}
$$
both coefficients are non-zero whenever $\v n\ge\v0$ and $\v n-\v
l\ge\v 0$. Hence the space of solutions of the form
(\ref{reducedsolution}) has dimension at most one. On the other
hand we do have such a solution, namely $\Psi_i$ where $i$ is
chosen such that the apexpoint $\beta_i$ is equal to the
apexpoint $\beta-\psi(\v k)$.

\qed
\medskip

We now consider polynomial mod $p$ solutions for A-hypergeometric systems with
parameters $\alpha\in\bbbq^r$.

\begin{proposition}\label{modpSolutions2}
Let $\alpha\in\bbbq^r$ and let $D$ the common denominator of the
coordinates of $\alpha$. Let $p$ be a prime not dividing $D$. Let
$\rho\is -p^{-1}\mod{D}$ if $D>1$ and $\rho=1$ if $D=1$. Let $s$
be the signature of $A$ and $\rho\alpha$. Suppose that the A-hypergeometric system
we consider is irreducible. Then, when $p$ is sufficiently large,
the polynomial mod $p$ solutions of the A-hypergeometric system with parameters
$A,\alpha$ is a free $\bbbf_p[v_i^{p}]$-module of rank $s$.
\end{proposition}

{\bf Proof}. Let $\v k=(1+p\rho)\alpha$. Notice that $\v
k\in\bbbz^r$ and $\v k\is \alpha\mod{p}$. So it suffices to look
at the mod $p$ A-hypergeometric system with parameters $A,\v k$. In Proposition
\ref{modpSolutions1} we saw that these solutions form a free
module of rank $s'$ where $s'$ is the signature of $A$ and $\v
k/p$. Let $\delta$ be the minimal distance of the points of
$\rho\alpha+\bbbz^r$ to the faces of $C(A)$. Suppose $\delta=0$.
Then there is a point $\rho\alpha+\v k$ with $\v k\in\bbbz^r$
contained in a face of $C(A)$. Choose $\mu\in\bbbz$ such that
$\mu\rho\is1\mod{D}$. Then $\mu(\rho\alpha+\v k)=\alpha+\mu\v k+(\mu\rho-1)\alpha$ is
on a face of $C(A)$. This contradicts the irreducibility of our
A-hypergeometric system by Theorem \ref{irreducible}. So $\delta>0$. Let us
assume that $p$ is so large that $|\alpha/p|<\delta$. Then the
points of $(\rho\alpha+\bbbz^r)\cap C(A)$ and $(\v
k/p+\bbbz^r)\cap C(A)$ are in one-to-one correspondence given by
$\v x\sim\v y\iff |\v x-\v y|<\delta$. In particular the number of
apexpoints of both sets is equal, hence $s=s'$. This proves our
assertion.

\qed

\section{Proof of the main theorem}
This section is devoted to a proof of Theorem \ref{algebraic}.
Let notations be as in Theorem \ref{algebraic} and suppose we
consider an irreducible A-hypergeometric system with parameters $\alpha\in\bbbq^r$.
Let $p$ be a prime
which is large enough in the sense of Proposition \ref{modpSolutions2}.
Let $D$ be the common denominator of the elements of $\alpha$ and
$\rho\is-p^{-1}\mod{D}$. Then the statement that $\sigma(A,\rho\alpha)$
is maximal is equivalent to the statement that the A-hypergeometric system modulo $p$ has a
maximal $\bbbf(\v v^p)$-independent set of polynomial solutions. 

A fortiori the following two statements are equivalent:
\medskip

i) $\sigma(A,k\alpha)$ is maximal for every
$k$ with $1\le k<D$ and $\gcd(k,D)=1$\\
ii) modulo almost every prime $p$ the A-hypergeometric system modulo $p$ has a maximal set
of polynomial solutions modulo $p$.
\medskip

A famous conjecture, attributed to Grothendieck implies that statement (ii) is equivalent to
the following statement,
\medskip

iii) The A-hypergeometric system has a complete set of algebraic solutions.
\medskip

If Grothendieck's conjecture were proven we would be done here. Fortunately, in two papers
by N.M.Katz (\cite{12} and \cite{15}) Grothendieck's conjecture is proven in the case when
the system of differential equations is (a factor of) a Picard-Fuchs system, i.e. a system
of differential equations satisfied by the period integral on families of algebraic varieties.
More precisely we refer to Theorem 8.1(5) of \cite{15}, which states

\begin{theorem}[N.M.Katz, 1982] Suppose we have a system of partial linear
differential equations, as sketched above, whose $p$-curvature vanishes for
almost all $p$. Then, if the system is a subsystem of a Picard-Fuchs system,
the solution space consists of algebraic functions.
\end{theorem}

The above theorem is formulated in terms of vanishing $p$-curvature for almost all $p$,
but according to a Lemma by Cartier (Theorem 7.1 of \cite{15}) this is equivalent to the
system having a maximal set of independent polynomial solutions modulo $p$ for almost all $p$.

To finish the proof of Theorem \ref{algebraic} it remains to show that the A-hypergeometric
equations for $\alpha\in\bbbq^r$ do arise from algebraic geometry. We shall do so in
Sections \ref{pochhammer} and \ref{eulerintegral}, where we construct Euler type integrals
for the solutions of the A-hypergeometric
system. 

In an attempt to maintain the lowtech nature of this paper we finish this Section with a proof
of the (easier) implication (iii)$\Rightarrow$ (ii). Before doing so we need a few introductory
concepts from the theory of linear differential equations.

Let $k$ be a field which, in our case, is usually $\bbbq$ or $\bbbf_p$.
Consider the differential field $K=k(v_1,\ldots,v_N)=k(\v v)$ with
derivations $\partial_i={\partial\over\partial v_i}$ for $i=1,\ldots,N$.
The subfield $C_K\subset K$ of elements all of whose derivatives are zero,
is called the field of constants. When the characteristic of $k$ is zero
we have $C_K=k$, when the characteristic is $p>0$ we have $C_K=k(\v v^p)$.

Throughout this section we let ${\cal L}$ be a finite set
of linear partial differential operators
with coefficients in $K$, like the A-hypergeometric system operators when $k=\bbbq$.
Consider the differential ring $K[\partial_1,\ldots,
\partial_N]$ and let $({\cal L})$
be the left ideal generated by the differential operators of the system.
We assume that the quotient $K[\partial_i]/({\cal L})$ is a $K$-vector
space of finite dimension $d$. Throughout this section we also fix
a monomial $K$-basis $\partial^{\v b}=\partial_1^{b_1}\cdots\partial_N^{b_N}$
with $\v b\in B$ and where $B$ is a finite set of $N$-tuples in $\bbbz_{\ge0}^N$
of cardinality $d$. 

\begin{proposition}\label{wronskian} Let ${\cal K}$ be some differential extension of
$K$ with field of constants $C_K$. Let $f_1,\ldots,f_m\in{\cal K}$ be a set
of $C_K$-linear independent solutions of the system $L(f)=0,\ L\in{\cal L}$.
Then $m\le d$. Moreover, if $m=d$ the determinant 
$$W_B(f_1,\ldots,f_d)=\det(\partial^{\v b}f_i)_{\v b\in B;i=1,\ldots,d}$$
is nonzero. 
\end{proposition}

In case we have $d$ independent solutions we call $W_B$ the {\it Wronskian
matrix} with respect to $B$ and $f_1,\ldots,f_d$. Obviously, if $g_1,\ldots,g_d$
are $C_K$-linear dependent solutions then $W_B(g_1,\ldots,g_d)=0$.
\medskip

{\bf Proof}. Suppose that either $m>d$ or $m=d$ and $W_B=0$. In both cases
there exists a ${\cal K}$-linear relation between the vectors
$df_i:=(\partial^{\v b}f_i)_{\v b\in B}$ for $i=1,2,\ldots,m$. 
Choose $\mu<m$ maximal such that $df_i,\ i=1,\ldots,\mu$ are ${\cal K}$-linear
independent. Then, up to a factor, the vectors $df_i,\ i=1,\ldots,\mu+1$
satisfy a unique dependence relation $\sum_{i=1}^{\mu+1}A_idf_i=0$
with $A_i\in{\cal K}$ not all zero. For any $j$ we can apply the operator
$\partial_j$ to this relation to obtain 
$$\sum_{i=1}^{\mu+1}\partial_j(A_i)df_i+A_i\partial_j(df_i)=0.$$
Since $\partial_j\partial^{\v b}$ is a $K$-linear combination of
the elements $\partial^{\v b},\v b\in B$ in $K[\partial_i]/({\cal L})$
there exists a $d\times d$-matrix $M_j$ with elements in $K$
such that $\partial_j(df_i)=df_i\cdot M_j$. Consequently 
$\sum_{i=1}^{\mu+1}A_i\partial_j(df_i)=\sum_{i=1}^{\mu+1}A_if_i\cdot M_j=0$
and so we are left with
$\sum_{i=1}^{\mu+1}\partial_j(A_i)df_i=0$. Since the relation between
$df_i,i=1,\ldots,\mu+1$ is unique up to factor there exists $\lambda_j\in K$
such that $\partial_j(A_i)=\lambda_jA_i$ for all $i$. Suppose $A_1\ne0$.
Then this implies that $\partial_j(A_i/A_1)=0$ for all $i$ and all $j$.
We conclude that $A_i/A_1\in C_K$ for all $i$. Hence there is a relation
between the $df_i$ with coefficients in $C_K$. A fortiori there is a 
$C_K$-linear relation between the $f_i$. This contradicts our assumption
of independence of $f_1,\ldots,f_m$. 

So we conclude that $m\le d$ and if $m=d$ then $W_B\ne0$.

\qed
\medskip

\begin{proposition}\label{algebraic2modP}
Suppose the system of equations $L(y)=0,\ L\in{\cal L}$ has only algebraic
solutions and that they form a vector space of dimension $d$.
Then for almost all $p$ the system of equations modulo $p$
has a $\bbbf(\v v^p)$-basis of $d$ polynomial solutions in $\bbbf(\v v)$.
\end{proposition}

{\bf Proof}. Let $f_1,\ldots,f_d$ be a basis of algebraic solutions. 
Choose a point $\v q\in\bbbq^N$ such that $f_i$ are all analytic near
the point $\v q$. Then $f_1,\ldots,f_d$ can be considered as power series
expansions in $\v v-\v q$. According to Eisenstein's theorem for powerseries
of algebraic functions we have that the coefficients of the $f_i$ can be
reduced modulo $p$ for almost all $p$. Moreover, let $\partial^{\v b},\v b\in B$
be a monomial basis of $K[\partial_i]/({\cal L})$. Then the Wronskian
determinant $W_B(f_1,\ldots,f_d)$ is non-zero. So for almost all $p$ the
powerseries $f_i$ can be reduced modulo $p$ and moreover, $W_B(f_1,\ldots,f_d)
\not\is0\mod{p}$. Hence, for almost all $p$ the powerseries $f_i\mod{p}$
are linearly independent over the quotient field of $\bbbf[[(\v v-\v q)^p]]$,
the power series in $(\v v-\v q)^p$. 

Fix one such prime $p$.
Let $P$ be the set $\{(b_1,\ldots,b_N)\in\bbbz^N\ |\ 0\le b_i<p\ {\rm for}
\ i=1,\ldots,N\}$. Every solution $f$ can be written in the form
$$f\is\sum_{\v b\in P}a_{\v b}(\v v-\v q)^{\v b}\mod{p},$$
where $a_{\v b}\in \bbbf[[(\v v-\v q)^p]]$. For every $L\in{\cal L}$
we have that
$$\sum_{\v b\in P}a_{\v b}L(\v v-\v q)^{\v b}\is0\pmod{p}.$$
Let $Q$ be the quotient field of $\bbbf[[(\v v-\v q)^p]]$.
The $Q$-linear relations between the polynomials
$L(\v v-\v q)^{\v b}$ for every $L$ form a vector space of dimension $d$
since the space of solutions mod $p$ has this dimension. Moreover
the space of $Q$-linear relations between the polynomials $L(\v v-\v q)^{\v b}$
is generated by $\bbbf((\v v-\v q)^p)$-linear
relations or, what amounts to the same, $\bbbfa(\v v^p)$-linear relations.

\qed

\section{Pochhammer cycles}\label{pochhammer}
In the construction of Euler integrals one often uses so-called twisted homology cycles.
In \cite{4} this is done on an abstract level, in \cite{19} it is done more explicitly.
In this paper we prefer to follow a more concrete approach by constructing a closed cycle of
integration such that the (multivalued) integrand can be chosen in a continuous manner
and the resulting integral is non-zero. For the ordinary Euler-Gauss function this is
realised by integration over the so-called Pochhammer contour. Here we construct its
$n$-dimensional generalisation. In Section \ref{eulerintegral} we use it to define an
Euler integral for A-hypergeometric functions.

Consider the hyperplane $H$ given by $t_0+t_1+\cdots+t_n=1$ in $\bbbc^{n+1}$ and the affine
subspaces $H_i$ given by $t_i=0$ for $(i=0,1,2,\ldots,n)$. Let $H^o$ be the complement in
$H$ of all $H_i$. We construct an $n$-dimensional
real cycle $P_n$ in $H^o$ which is a generalisation of the ordinary 1-dimensional Pochhammer
cycle (the case $n=1$). When $n>1$ it has the property that its homotopy class in $H^o$ is non-trivial,
but that its fundamental group is 
trivial. One can find a sketchy discussion of such cycles in \cite[Section 3.5]{16}.

Let $\epsilon$ be a positive but sufficiently small real number.
We start with a polytope $F$ in $\bbbr^{n+1}$ given by the inequalities
$$|x_{i_1}|+|x_{i_2}|+\cdots+|x_{i_k}|\le 1-(n+1-k)\epsilon$$
for all $k=1,\ldots,n+1$ and all $0\le i_1<i_2<\cdots<i_k\le n$. Geometrically this
is an $n+1$-dimensional octahedron with the faces of codimension $\ge2$ sheared off.
For example in the case $n=2$ it looks like

\centerline{\includegraphics[height=6cm]{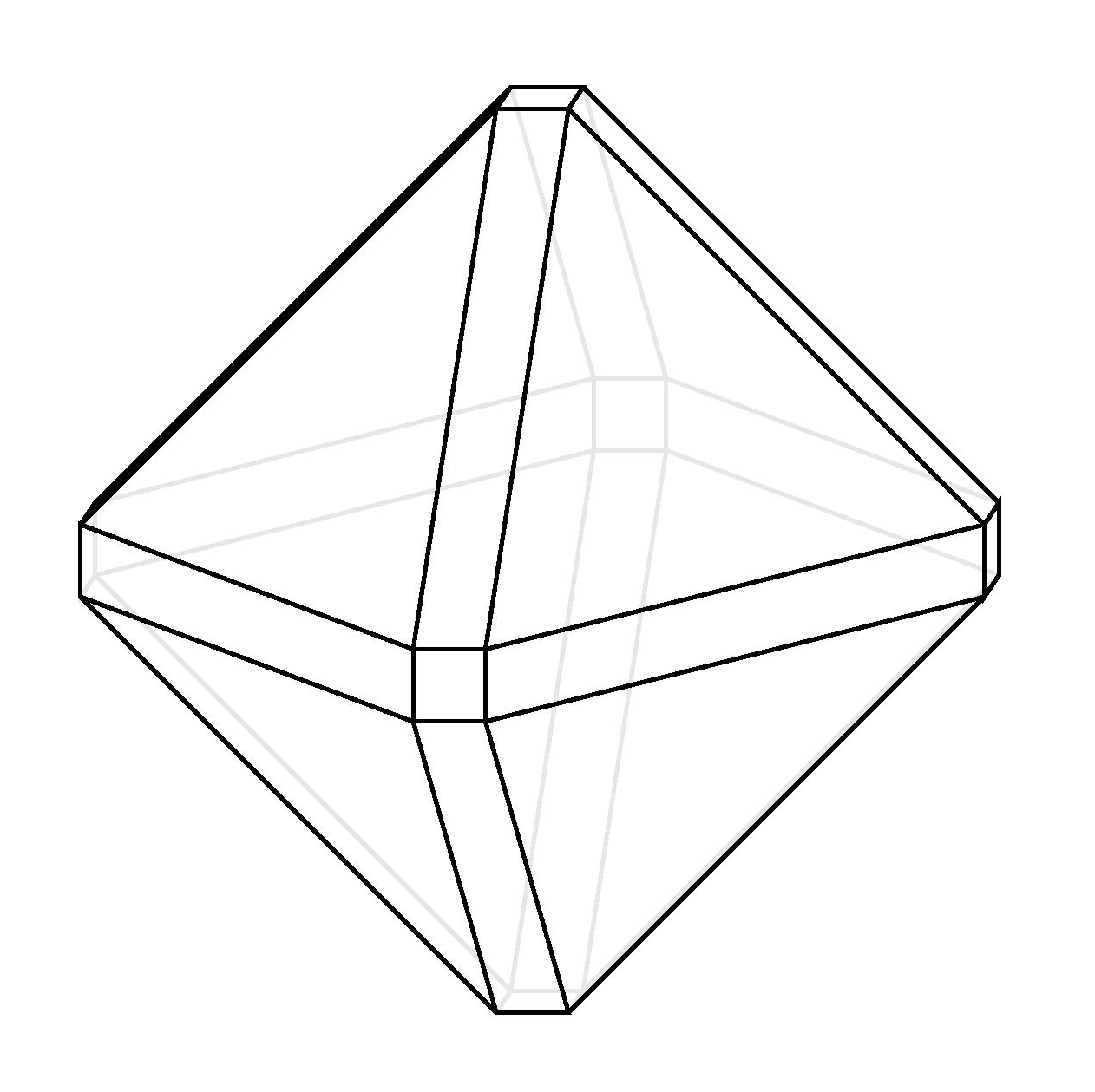}}

The faces of $F$ can be enumerated by vectors $\mu=(\mu_0,\mu_1,\ldots,\mu_n)\in\{0,\pm1\}^{n+1}$,
not all $\mu_i$ equal to $0$,
as follows. Denote $|\mu|=\sum_{i=0}^n|\mu_i|$. The face corresponding to $\mu$ is defined by
\begin{eqnarray*}
F_{\mu}&:& \mu_0 x_0+\mu_1 x_1+\cdots+\mu_n x_n=1-(n+1-|\mu|)\epsilon, \quad \mu_j x_j\ge\epsilon
\ {\rm whenever}\ \mu_j\ne0\\
&&|x_j|\le \epsilon\ {\rm whenever}\ \mu_j=0.
\end{eqnarray*}
Notice that as a polytope $F_{\mu}$ is isomorphic to $\Delta_{|\mu|-1}\times I^{n+1-|\mu|}$
where $\Delta_r$ is the standard $r$-dimensional simplex and $I$ the unit real interval. 
Notice in particular that we have $3^n-1$ faces. 

The $n-1$-dimensional side-cells of $F_{\mu}$ are easily described. Choose an index $j$
with $0\le j\le n$. If $\mu_j\ne0$ we set $\mu_j x_j=\epsilon$, if $\mu_j=0$ we set
either $x_j=\epsilon$ or $x_j=-\epsilon$. As a corollary we see that two faces $F_{\mu}$ and
$F_{\mu'}$ meet in an $n-1$-cell if and only if there exists an index $j$ such that $|\mu_j|\ne|\mu'_j|$
and $\mu_i=\mu'_i$ for all $i\ne j$. 

The vertices of $F$ are the points with one coordinate equal to $\pm(1-n\epsilon)$ and all
other coordinates $\pm\epsilon$.

We now define a continuous and piecewise smooth map $P:\cup_{\mu}F_{\mu}\to H$ as follows. Suppose the
point $(x_0,x_1,\ldots,x_n)$ is in $F_{\mu}$. Then its image under $P$ is defined as
\begin{equation}\label{intermediate}
{1\over y_0+y_1+\cdots+y_n}(y_0,y_1,\ldots,y_n)
\end{equation}
where $y_j=\mu_j x_j$ if $\mu_j\ne0$ and $y_j=E_{\epsilon}(x_j)$ if $\mu_j=0$. Here 
$E_{\epsilon}(x)=\epsilon e^{\pi i (1-x/\epsilon)}$. 
When $\epsilon$ is sufficiently small we easily check that $P$ is injective. We define
our $n$-dimensional Pochhammer cycle $P_n$ to be its image. 

\begin{proposition}\label{generalbeta}
Let $\beta_0,\beta_1,\ldots,\beta_n$ be complex numbers.
Consider the integral
$$B(\beta_0,\beta_1,\ldots,\beta_n)=
\int_{P_n}\omega(\beta_0,\ldots,\beta_n)$$
where 
$$\omega(\beta_0,\ldots,\beta_n)=t_0^{\beta_0-1}t_1^{\beta_1-1}\cdots t_n^{\beta_n-1}
\ dt_1\wedge dt_2 \wedge\cdots\wedge dt_n.$$
Then, for a suitable choice of the multivalued integrand, we have
$$B(\beta_0,\ldots,\beta_n)={1\over \Gamma(\beta_0+\beta_1+\cdots+\beta_n)}
\prod_{j=0}^n(1-e^{-2\pi i\beta_j})\Gamma(\beta_j).$$
\end{proposition}
{\bf Proof} 
The problem with $\omega$ is its multivaluedness. This is precisely the reason for
constructing the Pochhammer cycle $P_n$. Now that we have our cycle we solve the problem by making
a choice for the pulled back differential form $P^*\omega$ and integrating it over
$\partial F$.  Furthermore, the integral will not depend on the choice of $\epsilon$.
Therefore we let $\epsilon\to0$. In doing so we assume that the real parts of all $\beta_i$
are positive. The Proposition then follows by analytic continuation of the
$\beta_j$. 

On the face $F_{\mu}$ we define $T:F_{\mu}\to\bbbc$ by
$$T:(x_0,x_1,\ldots,x_n)=\prod_{\mu_j\ne0} |x_j|^{\beta_j-1}e^{\pi i(\mu_j-1)\beta_j}\prod_{\mu_k=0}
\epsilon^{\beta_j-1} e^{\pi i(x_j/\epsilon-1)(\beta_j-1)}.$$
This gives us a continuous function on $\partial F$. For real positive $\lambda$ we define
the complex power $\lambda^z$ by $\exp(z\log\lambda)$. 
With the notations as in (\ref{intermediate}) we have $t_i=y_i/(y_0+\cdots+y_n)$ and,
as a result,
$$dt_1\wedge dt_2 \wedge\cdots\wedge dt_n
=\sum_{j=0}^n (-1)^j y_jdy_0\wedge\cdots\wedge \check{dy_j}\wedge\cdots dy_n$$
where $\check{dy_j}$ denotes suppression of $dy_j$. 
It is straightforward to see that integration of $T(x_0,\ldots,x_n)$ over $F_{\mu}$ with $|\mu|<n+1$
gives us an integral of order $O(\epsilon^{\beta})$ where $\beta$ is the minimum
of the real parts of all $\beta_j$. Hence they tend to $0$ as $\epsilon\to0$. It
remains to consider the cases $|\mu|=n+1$. Notice that $T$ restricted to such
an $F_{\mu}$ has the form
$$T(x_0,\ldots,x_n)=
\prod_{j=0}^n e^{\pi i(\mu_j-1)\beta_j}|x_j|^{\beta_j-1}.$$
Furthermore, restricted to $F_{\mu}$ we have
$$\sum_{j=0}^n (-1)^j y_jdy_0\wedge\cdots\wedge \check{dy_j}\wedge\cdots dy_n
=dy_1\wedge dy_2 \wedge\cdots\wedge dy_n$$
and $y_0+y_1+\cdots+y_n=1$. Our integral over $F_{\mu}$ now reads
$$\prod_{j=0}^n\mu_j e^{\pi i(\mu_j-1)\beta_j}\int_{\Delta}(1-y_1-\ldots-y_n)^{\beta_0-1}y_1^{\beta_1-1}
\cdots y_n^{\beta_n-1}dy_1\wedge \cdots\wedge dy_n$$
where $\Delta$ is the domain given by the inequalities $y_i\ge\epsilon$ for
$i=1,2,\ldots,n$ and $y_1+\cdots+y_n\le 1-\epsilon$. The extra factor $\prod_j\mu_j$
accounts for the orientation of the integration domains. The latter integral is a
generalisation of the Euler beta-function integral. Its value is
$\Gamma(\beta_0)\cdots\Gamma(\beta_n)/\Gamma(\beta_0+\cdots+\beta_n)$. 
Adding these evaluation over all $F_{\mu}$ gives us our assertion.
\qed
\medskip

For the next section we notice that if $\beta_0=0$ the subfactor $(1-e^{-2\pi i\beta_0})
\Gamma(\beta_0)$ becomes $2\pi i$.

\section{An Euler integral for A-hypergeometric functions}\label{eulerintegral}
We now adopt the usual notation from A-hypergeometric functions. Define 
$$I(A,\alpha,v_1,\ldots,v_N)=\int_{\Gamma}{\v t^{\alpha}\over 1-\sum_{i=1}^Nv_i\v t^{\v a_i}}
\ {dt_1\over t_1}\wedge{dt_2\over t_2}\wedge\cdots\wedge{dt_r\over t_r},$$
where $\Gamma$ is an $r$-cycle which
doesn't intersect the hyperplane $1-\sum_{i=1}^Nv_i\v t^{\v a_i}=0$ for an open subset of
$\v v\in\bbbc^N$ and such that
the multivalued integrand can be defined on $\Gamma$ continuously and such that
the integral is not identically zero. We shall specify
$\Gamma$ in the course of this section. 

First note that an integral such as this satisfies the A-hypergeometric equations easily. The
substitution $t_i\to \lambda_i t_i$ shows that
$$I(A,\alpha,\lambda^{\v a_1}v_1,\ldots,\lambda^{\v a_n}v_N)=\lambda^{\alpha}I(A,\alpha,
v_1,\ldots,v_N).$$
This accounts for the homogeneity equations. For the "box"-equations, write $\v l\in L$
as $\v u-\v w$ where $\v u,\v w\in\bbbz_{\ge0}^N$ have disjoint supports.
Then

$$\Box_{\v l}I(A,\alpha,\v v)=|\v u|!\int_{\Gamma}{\v t^{\alpha+\v \sum_i u_i\v a_i}-
t^{\alpha+\v \sum_i w_i\v a_i}\over (1-\sum_{i=1}^Nv_i\v t^{\v a_i})^{|\v u|+1}}
\ {dt_1\over t_1}\wedge{dt_2\over t_2}\wedge\cdots\wedge{dt_r\over t_r}$$

where $|\v u|$ is the sum of the coordinates of $\v u$, which is equal to $|\v w|$
since $|\v u|-|\v w|=|\v l|=\sum_{i=1}^Nl_ih(\v a_i)=h(\sum_i l_i\v a_i)=0$.
Notice that the numerator in the last integrand vanishes because $\sum_iu_i\v a_i
=\sum_iw_i\v a_i$. So $\Box_{\v l}I(A,\alpha,\v v)$ vanishes.

We now specify our cycle of integration $\Gamma$. Choose $r$ vectors in $A$ such that
their determinent is $1$. After permutation of indices
and change of coordinates if necessary we can assume that $\v a_i=\v e_i$ for $i=1,\ldots,r$
(the standard basis of $\bbbr^r$). Our integral now acquires the form
$$\int_{\Gamma}{\v t^{\alpha}\over 1-v_1t_1-\cdots-v_rt_r-\sum_{i=r+1}^Nv_i\v t^{\v a_i}}
\ {dt_1\over t_1}\wedge{dt_2\over t_2}\wedge\cdots\wedge{dt_r\over t_r}.$$
Perform the change of variables $t_i\to t_i/v_i$ for $i=1,\ldots,r$. Up to a factor $v_1^{\alpha_1}
\cdots v_r^{\alpha_r}$ we get the integral
$$\int_{\Gamma}{\v t^{\alpha}\over 1-t_1-\cdots-t_r-\sum_{i=r+1}^Nu_i\v t^{\v a_i}}
\ {dt_1\over t_1}\wedge{dt_2\over t_2}\wedge\cdots\wedge{dt_r\over t_r},$$
where the $u_i$ are Laurent monomials in $v_1,\ldots,v_N$. Without loss of generality
we might as well assume that $v_1=\ldots=v_r=1$ so that we get the integral
$$\int_{\Gamma}{\v t^{\alpha}\over 1-t_1-\cdots-t_r-\sum_{i=r+1}^Nv_i\v t^{\v a_i}}
\ {dt_1\over t_1}\wedge{dt_2\over t_2}\wedge\cdots\wedge{dt_r\over t_r}.$$
For the $r$-cycle $\Gamma$ we choose the projection of the Pochhammer $r$-cycle
on $t_0+t_1+\cdots+t_r=1$ to $t_1,\ldots,t_r$ space. Denote it by $\Gamma_r$.
By keeping the $v_i$ sufficiently
small the hypersurface $1-t_1-\cdots-t_r-\sum_{i=r+1}^Nv_i\v t^{\v a_i}=0$ does not
intersect $\Gamma_r$. 

To show that we get a non-zero integral we set $\v v=\v 0$ and use the evaluation
in Proposition \ref{generalbeta}. We see that it is non-zero if all $\alpha_i$ have
non-integral values. When one of the $\alpha_i$ is integral we need to proceed with
more care. 

We develop the integrand in a geometric series and integrate it over $\Gamma_r$.
We have
\begin{eqnarray*}
&&{\v t^{\alpha}\over 1-t_1-\cdots-t_r-\sum_{i=r+1}^Nv_i\v t^{\v a_i}}\\
&&=\sum_{m_{r+1},\ldots,m_N\ge0}{|m|\choose m_{r+1},\ldots,m_N}
{\v t^{\alpha+m_{r+1}\v a_{r+1}+\cdots+m_N\v a_N}\over (1-t_1-\cdots-t_r)^{|m|+1}}
\ v_{r+1}^{m_{r+1}}\cdots v_N^{m_N}
\end{eqnarray*}

where $|m|=m_{r+1}+\cdots+m_N$.
We now integrate over $\Gamma_r$ term by term. For this we use Proposition
\ref{generalbeta}. We infer that all terms are zero if and only if there exists $i$
such that the $i$-th coordinate of $\alpha$ is integral and positive and the $i$-th
coordinate of each of $\v a_{r+1},\ldots,\v a_N$ is non-negative. In particular this means that
the cone $C(A)$ is contained in the halfspace $x_i\ge0$. Moreover, the points
$\v a_j=\v e_j$ with $j\ne i$ and $1\le j\le r$ are contained in the subspace $x_i=0$, so they span
(part of) a face of $C(A)$. The set $\alpha+\bbbz^r$ has non-trivial intersection 
with this face because $\alpha_i\in\bbbz$. From Theorem \ref{irreducible} it follows
that our system is reducible, contradicting our assumption of irreducibility.

So in all cases we have that the Euler integral is non-trivial. By irreducibility of
the A-hypergeometric system all solutions of the hypergeometric system can be given
by linear combinations of period integrals of the type $I(A,\alpha,\v v)$ (but with
different integration cycles).

\noindent
Department of Mathematics\\
Universiteit Utrecht\\
P.O. Box 80010, NL-3508 TA\\
Utrecht, The Netherlands\\
email: f.beukers@uu.nl 
\end{document}